\documentclass[
	a4paper, 
	10pt, 
	twoside, 
]{LTJournalArticle}
\usepackage{lipsum}
\DeclareOldFontCommand{\bf}{\normalfont\bfseries}{\mathbf}

\usepackage{amsmath}
\usepackage{algorithm}
\usepackage{algorithmic}
\usepackage{subcaption}
\usepackage{tikz}
\usepackage{graphicx}
\usepackage{graphbox} 
\usepackage{array}
\usetikzlibrary{shapes,patterns,arrows}
\usetikzlibrary{decorations.pathreplacing}

\newcommand{\mynorm}[1]{\| #1 \|}
\runninghead{} 

\footertext{}

\title{Dynamic image reconstruction in MPI with RESESOP-Kaczmarz} 

\author{%
	Marius Nitzsche, Bernadette N Hahn\\
}

\date{\footnotesize\textsuperscript{Department of Mathematics, University of Stuttgart, D-70569 Stuttgart, Germany}}

\renewcommand{\maketitlehookd}{%
\begin{abstract}
In Magnetic Particle Imaging (MPI), it is typically assumed that the studied specimen is stationary during the data acquisition. 
In practical applications however, the searched-for 3D distribution of the magnetic nanoparticles might show a dynamic behavior, caused by e.g. breathing or movement of the blood.
Neglecting those dynamics during the reconstruction step results in motion artifacts and a reduced image quality.

This article
addresses the challenge of capturing high quality images in the presence of motion.
A promising technique provides the Regularized Sequential Subspace Optimization (RESESOP) algorithm, which takes dynamics as model inexactness into account, significantly improving reconstruction compared to standard static algorithms like regularized Kaczmarz.
Notably, this algorithm operates with minimal prior information and the method allows for subframe reconstruction, making it suitable for scenarios with rapid particle movement. 
The performance of the proposed method 
is demonstrated on both simulated and real data sets.
\end{abstract}
}

\begin{document}

\maketitle 



\section{Introduction}
Magnetic Particle Imaging (MPI), initially proposed by Gleich and Weizenecker \cite{Gleich2005}, is an emerging imaging modality that enables insights into a specimen by measuring the response of superparamagnetic nanoparticles to an applied magnetic field.
To this end particles are injected with a certain tracer concentration $c: \Omega \mapsto \mathbb{R}_+ \cup {0}$ into the area of interest, which is called the field of view $\Omega \subset \mathbb{R}^3$.
Through layering two magnetic fields, i.e. the drive field and the selction field, a field-free point moves through the field of view along a Lissajous trajectory.
The time interval of one Lissajous trajectory is called \emph{frame}.
The dynamic field excitation caused by the moving field-free point results in a time-varying magnetization $M$ of the injected particles, which induces a measurable voltage $v^P$ over a time interval $I:=[0, \text{T}]$ in the receive coils.
To reconstruct $c$ it is necessary to separate it from $M$ by $M=c \Bar{m}$ with the particles' mean magnetic moment $\Bar{m}: \Omega \times I \mapsto \mathbb{R}^3$.
Knowing this, at the core of MPI lies the fundamental equation
\begin{align}\label{equ:basic}
    v^P(t) = &\int_{\Omega} c(x) \underbrace{(-a*\mu_0 p^R(x)^T \frac{\partial}{\partial t}\Bar{m}(x,t))}_{=:S(x,t)} \text{d}x, \; t\in I
\end{align}
where $\mu_0$ is the constant magnetic permeability of vacuum and $p^R: \Omega \mapsto \mathbb{R}^3$ is the receive coil sensitivity.
The analog filter $a : [-\text{T}, \text{T}] \to \mathbb{R}$ is applied to mitigate the direct feedthrough, a voltage induced by the applied magnetic field.

In MPI, two key challenges arise.
The determination of $\Bar{m}(x,t)$ from given tupels $(c,v^P)$ is called the \emph{calibration problem}.
The standard procedure relies on a time consuming measurement-based approach \cite{Knopp2012}.
Alternatively, model-based techniques are explored which offer a less resource-intensive path, although they may not deliver the same level of precision, see e.g. \cite{Albers2022}.

In this research, our primary focus is on the \emph{image reconstruction problem}.
In other words, the goal is to determine the tracer concentration $c(x), \, x\in \Omega$ from observed measurements $v^P(t), \, t\in I$ with a given system function.
This process necessitates the discretization of \eqref{equ:basic}, yielding a linear system
\begin{equation}\label{equ:dis}
    Ac=v, \quad A\in \mathbb{R}^{N\times M}, c\in \mathbb{R}^M, v \in \mathbb{R}^{N},
\end{equation}
with given system matrix $A$, measurement vector $v$ and searched-for discrete concentration $c$.
The dimensions are given by $N$, which is the number of considered temporal points multiplied by the number of receive coils being usually $3$, as well as the number of spatial points $M$ of the discretized field of view.

The most common solver in MPI is the Kaczmarz method with Tikhonov regularization, also called regularized Kaczmarz, which mainly consists of a fixed point iteration.
By iteratively projecting the approximate solution orthogonally onto affine subspaces, the method solves the linear system in \eqref{equ:basic}.
More insights can be found in the literature, see e.g. \cite{Knopp2012}.

\subsection{Dynamics in MPI}
The regularized Kaczmarz method works well if the searched-for concentration is stationary.
Be that as it may, many potential clinical applications rely on the reconstruction from dynamic MPI data, e.g. imaging blood flow \cite{Weizenecker2009}, tracking medical instruments \cite{Haegele2012} or monitoring of strokes \cite{Graeser2019, Ludewig2022}.
If the concentration is time-dependent and changes during the data acquisition, motion artefacts will arise in the reconstructions unless the dynamics is compensated for within the reconstruction step. This issue is exacerbated in the so-called multi-frame scenario, where the concentration is reconstructed from multiple data sets that are blockwise averaged over time to improve the signal-to-noise ratio and the spatial resolution of the image.

\begin{figure}[!htb]
    \centering
    \subcaptionbox{Phantom}
        {\includegraphics[width=.21\textwidth]{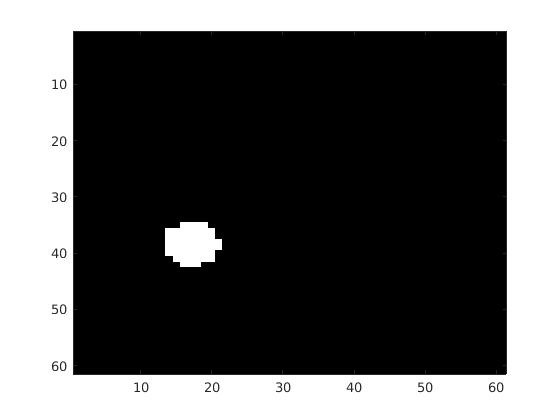}}
    \subcaptionbox{$\lambda = 8.5$: Averaging}
        {\includegraphics[width=.21\textwidth]{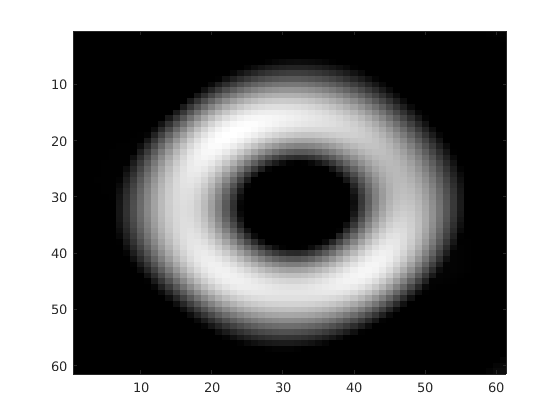}}

    \subcaptionbox{$\lambda = 2$: One frame}
        {\includegraphics[width=.21\textwidth]{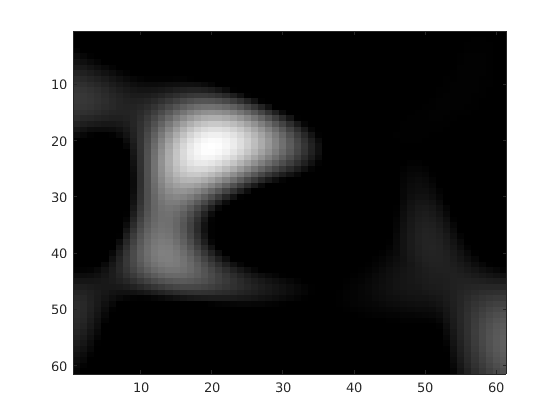}}
    \subcaptionbox{$\lambda = 8.5$: One frame}
        {\includegraphics[width=.21\textwidth]{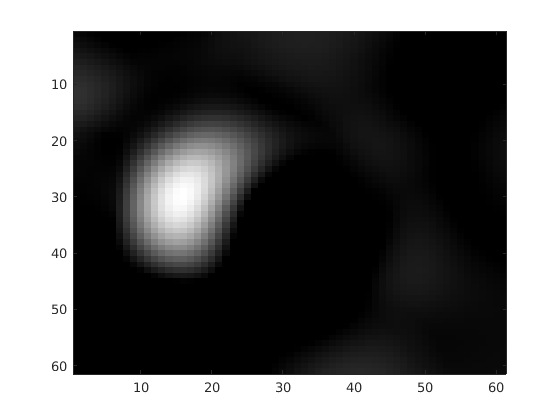}}
    \caption{Reconstruction results with regularized Kaczmarz from noise free data of a fast rotating cylinder.}
    \label{fig:regKacz}
\end{figure}

The regularized Kaczmarz method does not take dynamics of the concentration into account. 
\autoref{fig:regKacz} illustrates motion artefacts arising in single- as well as multi-frame reconstructions.
As example, we considered simulated noise free data for a fast rotating cylinder, more precisely it performs one full rotation during the acquisition of seven frames.
A 2D slice of the 3D ground truth at time instance $t=0$ is depicted in subfigure (a).
Reconstruction results from a single-frame data set with regularized Kaczmarz are shown in subfigures (c) and (d) for different regularization parameters $\lambda$ with the one in (d) being the optimal one.
In both cases, the computed solution is not sharp and includes motion artefacts. 
Even with the optimal parameter the method is not able to provide a good approximation to the ground truth while the averaging in the multi-frame scenario further increases the motion artifacts, see subfigure (b).

Consequently, it becomes imperative to account for dynamics, i.e. to address the image reconstruction task as a dynamic inverse problem. 
A universal, standardized regularization technique for dynamic inverse problems does not exist.
While a multitude of approaches have been proposed for individual traditional tomographic applications, only a few have been presented in the field of MPI in response to this challenge.

For instance, a strategy to account for periodic motion by pre-processing the measured data is developed in \cite{Gdaniec2017}.
By estimating the motion frequency associated with periodic motion, the authors effectively grouped data samples and thus generated static single-frame data sets for the individual states of the object.
Consequently, an image of each state can be reconstructed with standard static algorithms.
This approach has also been extended to multi-patch MPI \cite{Gdaniec2020}.

A dynamic forward model within MPI is offered in \cite{Brandt2021}, which is combined with a reconstruction method using temporal splines to increase temporal resolution in \cite{Brandt2022}.

Another approach is based on estimating the motion using optical flow, a computationally expensive method  that has been examined in applications like CT or MRI \cite{Burger2015}.
Its application in a dynamic MPI scenario was first proposed in \cite{kluth2019spatio} and in more detail in \cite{Brandt2023} using a stochastic primal-dual algorithm to jointly estimate the motion and reconstruct the image.

In this article we propose to apply the Regularized Sequential Subspace Optimization method (RESESOP) for inexact model operators \cite{Blanke2020} to the dynamic MPI problem.
With this method we take motion implicitly as model inexactness into account and are able to reconstruct images from motion corrupted data, which are of higher quality than those computed by regularized Kaczmarz, without the need of explicitly estimating the motion.
We will see that the only required a priori information can be gained from the measured MPI data itself. 
Furthermore, no special pre-processing of the data beyond the standard is required.

The article is organized as follows.
In section \ref{sec:methods} we first introduce the mathematical model of time-dependent concentration reconstruction in MPI and argue that dynamics induce model inexactness.
We then present the general concept of the RESESOP and transfer it to the dynamic reconstruction problem in MPI.
In particular, we discuss how the required a priori information on the model inexactness can be extracted directly from the measured MPI data.
Finally, we provide a detailed numerical evaluation of the method.
The respective test cases, including simulated as well as real data sets, are introduced in Section \ref{sec:data}.
In Section \ref{sec:results}, we present respective reconstruction results and analyze the performance of the proposed algorithm.

\section{Methods}\label{sec:methods}
As illustrated in the introduction, reconstruction algorithms need to account for dynamics of the studied concentration. In order to develop such algorithms, the time-dependency of the concentration first needs to be incorporated into the mathematical model of MPI, i.e. in equation \eqref{equ:basic} and \eqref{equ:dis}, respectively.

\subsection{The mathematical model of MPI with time-dependent concentrations}
Imaging a time-dependent particle concentration, i.e. $c:\Omega\times I \mapsto \mathbb{R}_+ \cup {0}$, implies that it is changing during the
measuring process. In the most extreme scenario, the time-scale of the concentration dynamics coincides
with the one of the data acquisition, i.e. 
\begin{equation}\label{equ:time_con}
    v^P(t) = \int_{\Omega} c(x,t) S(x,t) \text{d}x, \quad t\in I
\end{equation}
or in the discretized setting
\begin{equation}\label{equ:time_dis}
    A_t c_t=v_t \, \forall t \in [0,T], \quad A_t\in \mathbb{R}^{1\times M}, c_t\in \mathbb{R}^M, v_t \in \mathbb{R}^1.
\end{equation}
In other words, each single measure $v^P(t)$ (respectively $v_t$), $t$ being one
particular time instance, captures a different state of the concentration. Even with strong a priori assumptions a
dynamic concentration reconstruction on the same temporal granularity as defined by the sampling rate of the measurement process is hardly possible.

However, an important characteristic of MPI are the very fast measurement times which imply that the temporal resolution of the concentration is typically several orders coarser than the data sampling rate.
Thus, it is reasonable to assume that the searched-for concentration is piecewise constant in time.
This motivates to formally split the time interval $I$ into $N_\tau$ pairwise disjoint time subintervals $I_i:=[\tau_i, \tau_{i+1})$ with $i=0,1,\dots, N_\tau - 1$ as well as $\tau_0:=0$ and $\tau_{N_\tau}:=T$.
As a consequence it holds
\begin{align*}
    I=\dot\cup_{i=0,1,\dots,N_\tau - 1} I_i, \quad c(x,t)=c_{\tau_i} (x) \text{ for } t \in I_i.
\end{align*}
Both the time scale of the measurement and the coarser time scale of the dynamic concentration can then be coupled by a function
\begin{align*}
    \gamma: I \rightarrow \{ \tau_i, \ i=0,1,\dots, N_\tau - 1\},\  \gamma(t):= \tau_i \text{ if } t \in I_i.
\end{align*}
With this notation, the general dynamic reconstruction problem can be formulated as 
\begin{equation}\label{equ:time_con2}
    v^P(t) = \int_{\Omega} c(x,\gamma(t)) S(x,t) \text{d}x, \quad t\in I, 
\end{equation}
or in discretized form
\begin{equation}\label{equ:p_static}
    A_{i} c_{\tau_i}=v_{i} \quad \forall \ i\in \lbrace 0,1,\dots, N_\tau - 1 \rbrace
\end{equation}
with given $A_{i}\in \mathbb{R}^{(N_i\cdot 3)\times M}$, searched-for $c_{\tau_i} \in \mathbb{R}^M$ and measured $v_{i} \in \mathbb{R}^{(N_i \cdot 3)}$, where $N_i$ denotes the number of time samples per sub-interval $I_i$.

In Section \ref{Sec:II-II}, we will focus on solving the discretized dynamic inverse problem.
Considering each \emph{subproblem} $A_{i} c_{\tau_i}=v_i$ individually would formally allow the application of a classic static reconstruction method. 
However, in analogy to the multi-frame approach, the goal is to improve the SNR and the spatial resolution of the images.
In order to achieve this in the presence of motion, correlations between the states have to be exploited within the joint reconstruction step. 
Furthermore, if the length of a subinterval $I_i$ is smaller than the time required to measure one complete frame, the solution of the respective subproblem has to deal with a limited-data problem.  
The actual decomposition of $I$ into subintervals will depend on the specific application. Suitable choices are discussed in the following subsection.

\subsubsection{Choosing suitable subintervals}
In the multi-frame scenario, an intuitive idea is to let each subproblem coincide with one frame. Indeed, many concentration dynamics might be assumed stationary within
one frame. 
The actual time to acquire data for one frame is very small, e.g. for a 3D field of view, it takes 21.54 ms. Thus, for velocities lower than 10 centimeters per second, for instance, the travelled distance within the 21.54 ms is smaller than the diameter of a voxel size \cite{Kaul_Salamon_Knopp_Ittrich_Adam_Weller_Jung_2018}. 
Additionally, the usage of frames as subproblems simplifies the implementation within a MPI framework.
Real data are usually pre-processed and stored frame-by-frame. Thus, a dynamic reconstruction algorithm working on frames can be embedded directly in existing MPI software pipelines.  

However, if the examined structure changes faster than the acquisition of one frame, e.g. due to a pulsating motion or a fast rotation, this rigid coupling of time scales is not feasible and would still lead to motion artifacts.
A way to handle these scenarios is to choose smaller subproblems which are then underdetermined. 
Looking at it from a practical perspective, one timestep $t$ is equivalent to $0.8 \mu\text{s}$ in MPI.
Thus, choosing subproblems corresponding to either one frame, half a frame or a quarter of a frame seem adequate to represent and capture most relevant concentration dynamics.

\subsubsection{Interpreting the dynamics as model inexactness}

A common strategy in dynamic imaging problems is the incorporation of motion models to relate the different states of the object to each other,
see e.g [\cite{Hahn_2014},, 
\cite{Desbat}]. In our case, this means that instead of considering a series of independent concentrations $c_{\tau_i}$, it is assumed that they are all linked together by an underlying motion model $\Gamma$, for instance
\begin{align*}
    c(x,t)=c_0(\Gamma_t(x))
\end{align*}
with a reference configuration $c_0$, e.g. the concentration
at the inital time of the data acquisition. 
Incorporating such a motion model into the equation \eqref{equ:time_con2} and using an appropriate change of coordinates yields an inverse problem for the reference concentration $c_0$ with a forward operator depending on the motion information $\Gamma$, i.e. in the discretized setting 
\begin{equation}\label{equ:p_static2}
    A_{i,\Gamma} c_{0}=v_{i} \quad \forall \ i\in \lbrace 0,1,\dots, N_\tau - 1 \rbrace 
\end{equation}
with $A_{i,\Gamma}\in \mathbb{R}^{(N_i\cdot 3)\times M}$.

A detailed derivation in the general context of inverse problems can be found in \cite{Hahn_2014}.
Thus, extracting the time-dependent concentrations $c(\cdot,t)$ is equivalent to extracting the motion information $\Gamma$ and the reference concentration $c_0$.

In general, the exact motion $\Gamma$ will be unknown, i.e. in practice, one can only use an approximate forward operator to determine $c_0$, e.g. the static model $A_i$ if no further information on the dynamics are available.
Therefore, we treat the dynamics in the following as inexactness in the forward model and account for this inexactness while recovering $c_0$ from the measured dynamic data $v_i, \ i=0,1,\dots,N_{\tau}-1$.

\subsection{ RESESOP method for solving inverse problems with inexact forward operator}\label{Sec:II-II}

In this work we propose to solve the dynamic reconstruction problem in MPI \eqref{equ:p_static2} with the RESESOP-Kaczmarz method developed in \cite{Blanke2020} for solving general inverse problems with inexact forward operator.

\begin{figure}[!htb]
\centering
  \begin{tikzpicture}[xscale=0.7,yscale = 0.6]  
    \draw [very thick, red](0,0) -- (6,0);
    \draw (0,0.7) -- (6,0.7);
    \draw (0,-0.7) -- (6,-0.7);
    \draw[fill=gray!20,nearly transparent]  (0,-0.7) -- (6,-0.7) -- (6,0.7) -- (0,0.7) -- cycle;
    \draw [red](6,0) node[right=3pt] {$ H(u,\alpha) $};
    \draw (6,0.7) node[right=3pt] {\footnotesize{$ H(u,\alpha+\xi) $}};
    \draw (6,-0.7) node[right=3pt] {\footnotesize{$ H(u,\alpha-\xi) $}};
    \draw[<->](2,0.7) -- (2,0);
    \draw (2,0.35) node[right=1pt]{\footnotesize{$ \xi $}};
    \draw[decorate,decoration={brace,mirror}](-0.1,0.7) -- (-0.1,-0.7);
    \draw (-0.2,0) node[left=1pt]{$ H(u,\alpha, \xi) $};
  \end{tikzpicture}
  \vspace{-0.7em}
  \caption{Illustration of a stripe}
    \label{fig:stripe}
\end{figure}
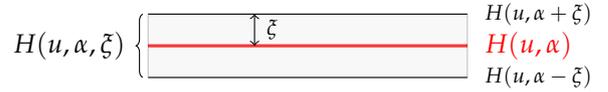

This method is based on the Sequential Subspace Optimization method (SESOP) \cite{Schoepfer2007}. It is an iterative method that approximates the searched-for solution $f$ of a given linear inverse problem $Bf=g$ by the following scheme: In the $n$-th step, compute   
$$ f_{(n+1)} = P_H(f_{(n)})$$
with the metric projection $P_H$ onto an intersection $H=\cap_{j \in J_n} H(u_j, \alpha_j)$ with hyperplanes
$$  H(u,\alpha)=\{ x \in X: \langle u,x \rangle = \alpha  \}.$$
In particular, each iteration has the following choices of freedom: 
\begin{itemize}
\item a finite index set $J_n$, denoting the number of search directions $u_j$ and thus the number of hyperplanes used in the $n$-th iteration,
\item search-directions $u_j$ and parameters $\alpha_j$ which are chosen in accordance to the right-hand side $g$ and the forward model $B$, for instance, in case $J=\lbrace n\rbrace$ a suitable choice is $u_n=B^*(Bf_{(n)}-g)$ as well as $\alpha_n=\langle Bf_{(n)}-g, g \rangle$. 
\end{itemize}
This method relies on knowing the exact forward model $B$ and on exactly known right-hand side $g$.
In practice however, only noisy data $g^\delta$ with noise level \mbox{$\mynorm{g-g^{\delta}}\leq \delta$} can be measured and, as in our case, only an inexact forward operator $B^\eta$ with inexactness level \mbox{$\mynorm{B-B^{\eta}}\leq \eta $} might be available for the reconstruction.
Thus, the idea of RESESOP (regularized SESOP) is to replace the hyperplanes $H(u,\alpha)$ by stripes $$H(u,\alpha, \xi):=\{ x \in X: |\langle u,x \rangle - \alpha| \leq \xi  \}, $$
whose width $\xi$ is chosen in dependence on the noise and inexactness levels $\delta$ and $\eta$, see \autoref{fig:stripe} for an illustration.
This idea was originally suggested in \cite{Schoepfer2009} for the noisy data case and was expanded to model imperfections in \cite{Blanke2020}.
The introduction of Morozov's discrepancy principle offers a stopping criterion; as a result, it can be proven under certain conditions on the choice of $J_n, u_i, \alpha_i$, that the method provides a regularized solution to the original inverse problem \cite{Schoepfer2009}.

When solving an inverse problem that shows a subproblem structure like (7), it is beneficial to include individual inexactness and noise levels $\eta_i, \, \delta_i$,
in our case
$$ \mynorm{v_{i}-v_{i}^{\delta}}\leq \delta_i , \quad \mynorm{A_{i}-A_{i}^{\eta}}\leq \eta_i, \quad i\in \lbrace 0,1,\dots,N_{\tau}-1\rbrace. $$
In particular, with such estimates, we obtain 
$$ \mynorm{A_{i}^{\eta} c - v^\delta_{i} }\leq \eta_i \, \rho + \delta_i$$
as upper bound for the residual on the set of solutions whose norm is bounded by $\rho>0$. 
For this case, Blanke et al proposed the so-called RESESOP-Kaczmarz method: Each iteration corresponds to a RESESOP iteration with subproblem specific inexactness and noise levels applied to every subproblem that does not yet satisfy the discrepancy principle. In this way, the algorithm combines all subproblems and incorporates information from all time intervalls $I_i$ simultaneously.

A general form of the RESESOP-Kaczmarz method is stated in algorithm \ref{alg:res}.
With $n$ being the current iteration index, the notation $[n]:= n \text{ mod }N_\tau$ identifies the corresponding subproblem \eqref{equ:p_static}. 
One full iteration consists of $N_\tau$ subiterations, each including a metric projection onto the intersection of stripes $H_{(n)}^{\eta,\delta}$.

\begin{algorithm}
\caption{RESESOP-Kaczmarz}
\label{alg:res}
\begin{algorithmic}
\STATE Choose $c_{(0)}$ and constant $\rho > 0$.\\
Let $n$ be the iteration index.
\WHILE {$c_{(n)} \neq c_{(n-N_\tau)}$ for $[n] = 0$}
\IF {$||A_{[n]}c_{(n)}-v_{[n]}^{\eta,\delta}|| \leq 1.001 \cdot (\eta_{[n]} \rho + \delta_{[n]})$}
    \STATE $c_{(n+1)}=c_{(n)}$
\ELSE
    \STATE Choose a finite index set $J_n \subset \{ 0,1,\dots,n\}$ 
    \STATE Choose search directions $u_{n,j}$ for all $j\in J_n$
    \STATE Define $H_{(n)}^{\eta,\delta}:=\underset{j\in J_n}{\cap} H (u_{n,j}, \alpha_{n,j}, \xi_{n,j})$
    \STATE with parameters $\alpha_{n,j}=\langle u_{n,j},c_{(n)} \rangle$ \\ and $\xi_{n,j}=(\eta_{[j]} \rho + \delta_{[j]})\mynorm{A_{[j]} c_{(n)} - v_{[j]}^{\eta,\delta}}$
    \STATE Compute $c_{(n+1)}=P_{H_{(n)}^{\eta,\delta}}(c_{(n)})$
\ENDIF
\ENDWHILE
\end{algorithmic}
\end{algorithm}

RESESOP-Kaczmarz does not require extensive or unrealistic a priori information to solve inverse problems with inexact forward model such as dynamic reconstruction problems. 
The estimates of the uncertainty levels are the only additional information used compared to standard static approaches. Consequently, when applying the method to dynamic problems, the motion information will be implicitly included in the inexactness levels. 

However, a sufficiently good estimate is crucial for the performance of RESESOP.
If the estimates are too large the convergence of the method can be fast but lead to a computed solution with a high error rate.
If the levels are estimated too low the noise and the model imperfections might not be compensated for sufficiently resulting in a noisy reconstruction with artefacts.
The extreme case, where the levels are set to zero, corresponds to a SESOP-Kaczmarz 
reconstruction without regularization and accounting for model inexactness, i.e. in particular to a static reconstruction where the Kaczmarz loop results in an averaging of the data.

Next, we discuss further adaptions of this general framework to the application in dynamic MPI, in particular how the required estimates can be obtained and how the search directions can be chosen.

\subsection{Application in MPI}
\subsubsection{Inexactness levels in dynamic MPI} \label{sec:in}
We first discuss how to determine estimates for the noise and inexactness levels. In a static multi-frame setting, the noise levels $\delta_i$ can
be determined by comparing the data vectors measured per frame.

The inexactness levels depend on the selected model for the reconstruction algorithm.
Due to the time consuming process of calibration in real MPI and since we assume no further prior information on the unknown motion $\Gamma$, we propose to use the static system matrix $A_i$ as inexact forward model. Without knowledge of $\Gamma$, it is difficult to compute a viable upper bound $\eta_i$ with $\mynorm{A_i-A_{i,\Gamma}}\leq \eta_i$. 

However, if we compare again data sets for different frames in a dynamic scenario, their difference captures noise as well as motion.
For this reason, we propose to jointly estimate both uncertainty levels together in form of levels $\zeta_i \in \mathbb{R}$ characterizing the total inexactness of the system, more precisely for each subproblem
\begin{equation*}
   \zeta_{i} \approx \eta_{i} \rho + \delta_{i}, \quad i\in \lbrace 0,1,\dots,N_{\tau}-1\rbrace.
\end{equation*}
These values are computed by comparing the data  
$v_{0}$ of the initial subproblem with data $v_{i}$ of the $i$-th subproblem.
If the size of each subproblem corresponds to one frame, it is simply $\zeta_{i}=\mynorm{v_0-v_i}$. 

However, in case of subframe-sized subproblems, this direct comparison is not possible but requires an additional interpolation step beforehand.
Each subproblem corresponds to a number of time points $t$ of a frame.
Consequently, the data fragment of the initial subproblem is comparable to the data acquired at the same time points of each frame, e.g. the uncertainty levels between the first subproblem of the first frame and the first subproblem of every other frame can be evaluated.
However, to compute an estimate for the uncertainty level of another subproblem of the same frame with the same method is impossible.
Therefore, we apply a cubic interpolation between comparable data fragments to acquire uncertainty levels for each subproblem.

An alternative method to estimate the uncertainty levels $\zeta_i$ is to make use of preliminary regularized Kaczmarz reconstructions.
To quantify the error between Kaczmarz's reconstructions, the mean-squared error MSE is used.
However, this approach depends on a priori reconstructions which will typically contain artefacts and thus can result in a low MSE. Thus, an extraction directly from measured data seems favorable.

\subsubsection{Choice of search directions}
The convergence properties of SESOP, and hence RESESOP-Kaczmarz, depend on the number and choice of the search-directions. 
On the one hand, a higher number of search directions results in more complex and computationally expensive iterations.
On the other hand, the algorithm requires in total less iterations to achieve the desired accuracy.
In this work we will focus on two search directions which is a compromise between computational expense for one iteration as well as total number of iterations.

More precisely, we choose in the $n$-th iteration step the finite index set $J_n=\{n_-,n\}$ with $n_-$ representing the last iteration in which the discrepancy principle did not yet hold.
The search directions $u_{n,j}$ with $j \in J_n$ are then chosen as
\begin{align*}
    u_{n,n_-}&= A_{{[n_-]}}^* (A_{{[n_-]}} c_{(n_-)} - v_{{[n_-]}}),\\
    u_{n,n}&= A_{{[n]}}^* (A_{{[n]}} c_{(n)} - v_{{[n]}}).
\end{align*}
In particular, this choice meets the criteria formulated in \cite{Blanke2020} that guarantee the convergence of RESESOP-Kaczmarz to a regularized solution of the underlying dynamic inverse problem.

\subsubsection{Positivity constraints}
In the MPI setting the reconstructed image approximates the concentration $c$, which is a vector consisting of real, positive numbers.
Therefore, a positivity constraint is a common technique to improve reconstruction algorithms \cite{Knopp2012}.
To that end, after each iteration every negative part is set to zero.

\section{Data}\label{sec:data}
Our test cases for the numerical evaluation comprise both real and simulated data, each offering valuable insights into the dynamic reconstruction problem in MPI and the performance of the proposed method.
This dual test set strategy allows for a robust evaluation under controlled conditions (simulated data) as well as real-world scenarios (real data).
\subsection{Simulated Data}\label{sec:datasim}
We first conduct a fully simulated numerical experiment.
This simulated dataset was generated using the model B3 presented by Kluth, Szwargulski, and Knopp \cite{KluthSzwargulskiKnopp2019}.

\begin{figure}[!htb]
\centering
\begin{subfigure}{.125\textwidth}
  \centering
  \includegraphics[width=\linewidth]{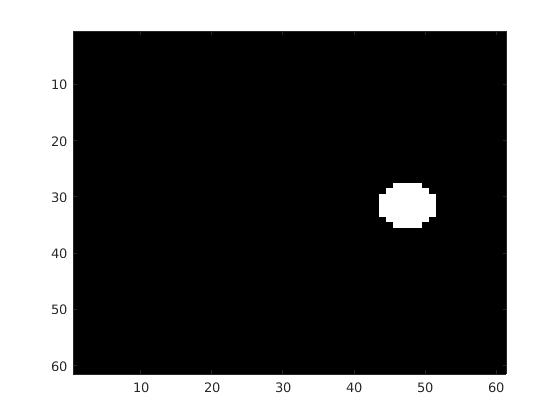}
  \caption{Frame 1}
\end{subfigure}%
\begin{subfigure}{.125\textwidth}
  \centering
  \includegraphics[width=\linewidth]{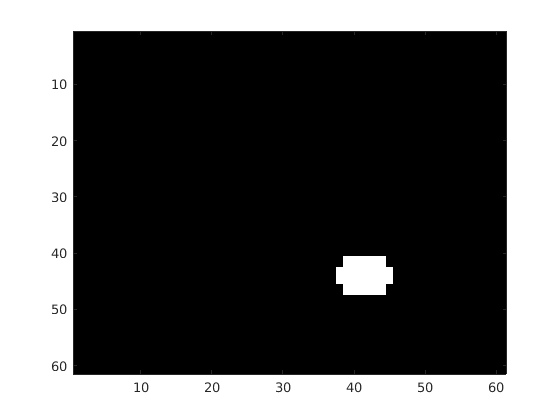}
  \caption{Frame 2}
\end{subfigure}%
\begin{subfigure}{.125\textwidth}
  \centering
  \includegraphics[width=\linewidth]{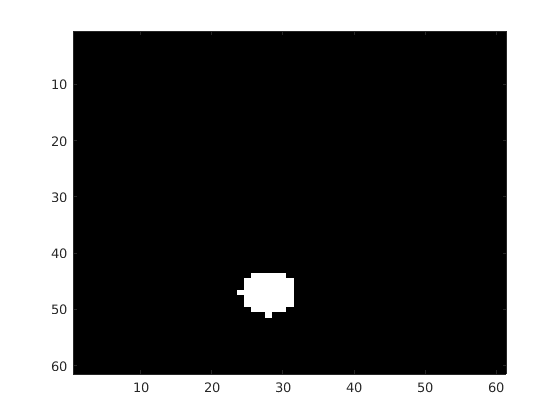}
  \caption{Frame 3}
\end{subfigure}%
\begin{subfigure}{0.125\textwidth}
  \centering
  \includegraphics[width=\linewidth]{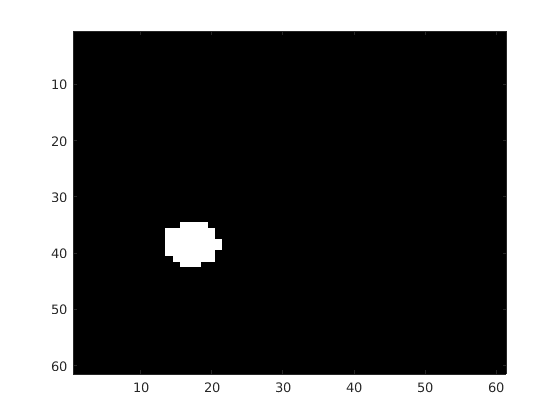}
    \caption{Frame 4}
\end{subfigure}\\
\begin{subfigure}{.125\textwidth}
  \centering
  \includegraphics[width=\linewidth]{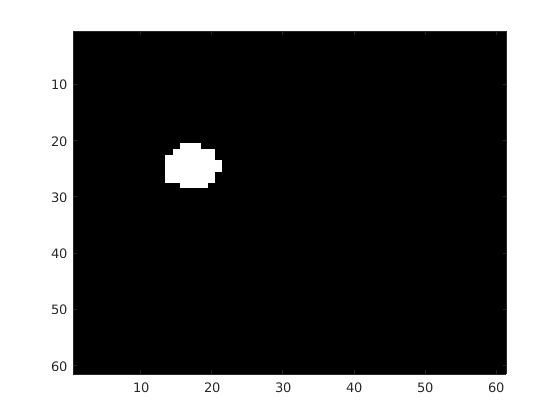}
  \caption{Frame 5}
\end{subfigure}%
\begin{subfigure}{.125\textwidth}
  \centering
  \includegraphics[width=\linewidth]{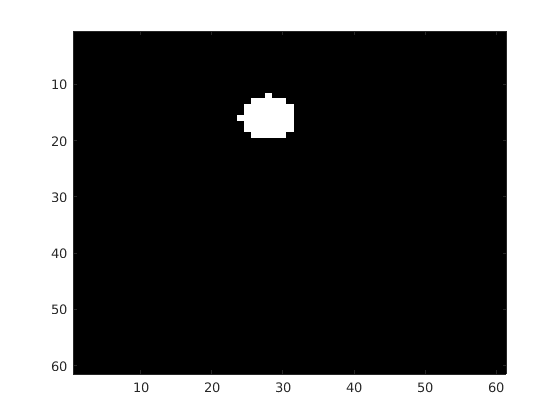}
  \caption{Frame 6}
\end{subfigure}%
\begin{subfigure}{.125\textwidth}
  \centering
  \includegraphics[width=\linewidth]{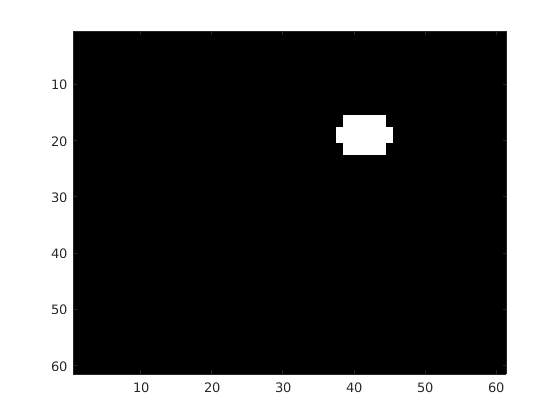}
  \caption{Frame 7}
\end{subfigure}%
\begin{subfigure}{0.125\textwidth}
  \centering
  \includegraphics[width=\linewidth]{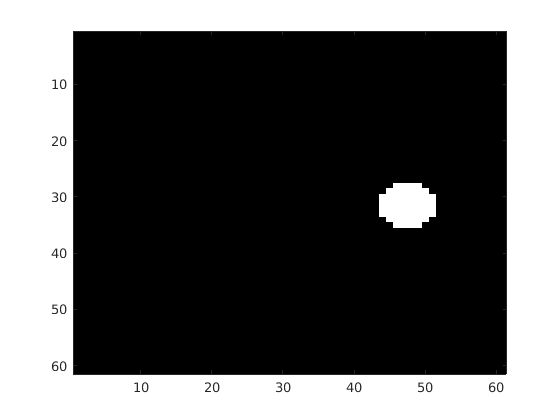}
    \caption{Frame 8}
\end{subfigure}\\
\caption{The course of the motion illustrated by the ground truth phantom at the beginning of each trajectory of the fast rotation (seven fpr).}
\label{fig:phantoms}
\end{figure}

The goal was to construct a similar setup as the available real data.
To that end, the phantom is a cylinder rotating through the field of view.
There are two different rotation speeds, depending on the number of frames per rotation (fpr).
They are equivalent to seven and 44 fpr.
\autoref{fig:phantoms} illustrates the underlying rotating motion of the object in the simulated test data case.

The simulated system matrix is derived from dynamic simulations of N\'{e}el-type particle magnetization dynamics.
Here, the uniaxial particle anisotropy varies with position within the field of view.
The anisotropy is minimal at the center but intensifies toward the boundaries.
This choice aligns with the assumptions in MPI, where the static component of the applied magnetic field is believed to induce physical rotation of particles and influence their combined anisotropy energy landscape.

The simulated data were generated on a $3$D grid of $61 \times 61 \times 5$ pixels, each with a resolution of $0.5$ mm.
A $2$D Lissajous excitation was used, with anisotropy constants reaching values as high as $1250 , \text{J/m}^3$.

To create the system matrix, we utilized a dedicated toolbox presented in Albers, Kluth, and Knopp \cite{AlbersKluthKnopp2022}.
To prevent inverse crime, the data were calculated using a slightly shifted, larger matrix, ensuring the integrity of our simulation experiments.
Furthermore, we added white Gaussian noise to the data corresponding to a signal-to-noise ratio of 10. 
\subsection{Real Data}\label{RealData}
The real data were 
provided by Brandt from the University of Hamburg and collected by a team at UKE Hamburg led by Knopp.
These data encode a glass capillary, which rotates at an average speed of seven Hz and was initially introduced in the study by Gdaniec et al. \cite{Gdaniec2017}.
This experiment represents a very fast motion.

Working with real data presents some additional challenges compared to simulated data. 
Typically, one lacks a ground truth which hinders a direct assessment of the accuracy of the reconstruction.
Moreover, real measurements usually have high degree of noise, adding complexity to the reconstruction task in the dynamic case where noise levels cannot be reduced by averaging over multiple frames.
This is further amplified in MPI due to the fact that the system function is measured as well.

These experiments were conducted using a pre-clinical MPI scanner from Bruker Biospin in Ettlingen.
The imaging setup involved 3D MPI measurements with a Lissajous excitation driven by three sinusoidal drive fields in the $x$-, $y$-, and $z$-directions.
These fields operated at frequencies of $f_x = \frac{2.5 \text{MHz}}{102}$, $f_y = \frac{2.5 \text{MHz}}{96}$, and $f_z = \frac{2.5 \text{MHz}}{99}$, with an amplitude of $14 \frac{\text{mT}}{\mu_0}$.
Each cycle of data acquisition took $21.54$ ms, and the induced signal was sampled at intervals of $0.8 \mu\text{s}$.
The measurements covered $15.625$ positions in a $25 \times 25 \times 25$ grid.
More detailed information can be found in \cite{Gdaniec2017}.

\subsection{Pre-processing}
To the real data, the standard pre-processing steps were applied.
After applying 
the Fourier transform, a frequency selection was performed, more precisely all frequencies below 80 k Hz and above 625 k Hz were eliminated.
This was followed by a SNR threshholding. 
Furthermore, a singular value decompostion as proposed in \cite{Kluth_2019} was performed to reduce the computational cost of reconstructing a three dimensional image.
With weigths computed from background measurements, the randomized singular value decomposition trimmed the system matrix to 5000 rows.

However, it is important to note that these conventional pre-processing steps are tailored to full-frame data. Thus, they may encounter limitations when applied to dynamic cases, where the dynamics requires a division of the data into sub-frames. 
In order to realize this, the raw data first have to be transformed into the time domain with an inverse Fourier transform.
Only then the frame-sized problems can be split into smaller subproblems. To this end, access to the original raw-data is crucial for the solution of reconstruction problems with fast dynamics. If only pre-processed data are provided, one is restricted to decompositions, where the time intervalls $I_i$ correlate to individual frames.

When working with simulated data, most of the pre-processing steps are not necessary.
Due to the 2D-Lissajous excitation, the system matrix is naturally smaller compared to the system matrix of real data.
Therefore, pre-processing steps to reduce computational work are not critical and some steps e.g. SNR threshholding not applicable to simulated data.

\section{Results and Discussions}\label{sec:results}
This section is dedicated to the numerical evaluation of the RESESOP-Kaczmarz method on both simulated and real dynamic MPI data. In particular, we investigate the influence of critical parameters of our algorithm on the result. 

All experiments are conducted with 3D data, but for clarity, each result is shown with a 2D slice.
In each experiment, the same frame and slice are depicted, ensuring a comparable view across all scenarios.

\subsection{Simulated Data}
We first present reconstruction results from simulated noisy data, cf. Section \ref{sec:datasim}, with two different rotation speeds - a very fast scenario of seven fpr as well as the scenario of 44 fpr.
Regarding the coupling of time scales, we consider here each frame as one subproblem, i.e. we are choosing the least intrusive version in the MPI pipeline.

To illustrate the motion compensation properties of the RESESOP-Kaczmarz method, we compare the results with those generated by the regularized Kaczmarz-algorithm which is the most commonly used algorithm within MPI. Since averaging data is not recommendable for dynamic data, all reconstructions with regularized Kaczmarz use one frame of data. 
As recommended in the literature, e.g. \cite{Knopp_2010}, the algorithm is stopped early after five iterations since the quality of the reconstruction deteriorates with more iterations. 
Furthermore, at each iteration positivity of the approximated solution is enforced.
The initial value was set to zero. The regularization parameter $\lambda$ is optimized by performing reconstructions for a wide range of parameters and choosing the one that leads to the best result.
\begin{figure}[!htb]
\centering
\begin{subfigure}{.23\textwidth}
  \centering
  \includegraphics[width=\linewidth]{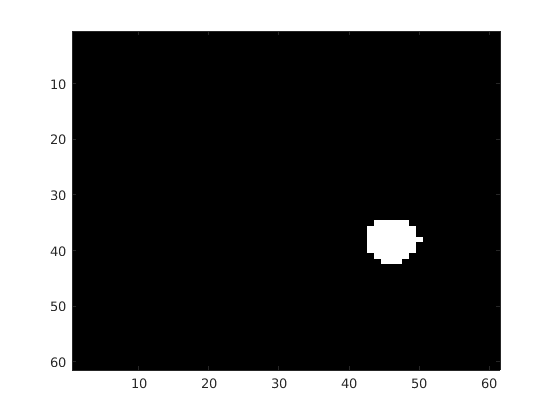}
  \caption{Ground truth}
\end{subfigure}%
\begin{subfigure}{0.23\textwidth}
  \centering
  \includegraphics[width=\linewidth]{figures/phantom_frame4_speed7.png}
    \caption{Ground truth}
    \label{fig:phantom_fast}
\end{subfigure}\\
\begin{subfigure}{.23\textwidth}
  \centering
  \includegraphics[width=\linewidth]{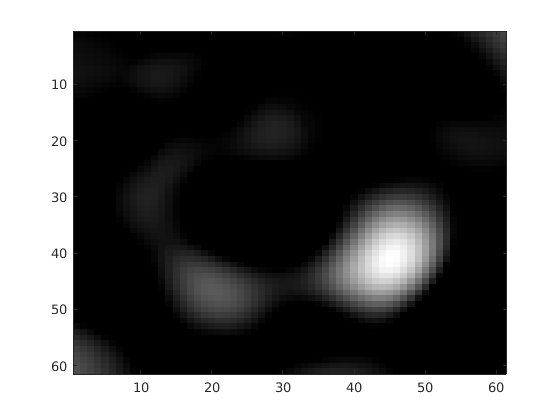}
  \caption{Regularized Kaczmarz}
\end{subfigure}%
\begin{subfigure}{0.23\textwidth}
  \centering
  \includegraphics[width=\linewidth]{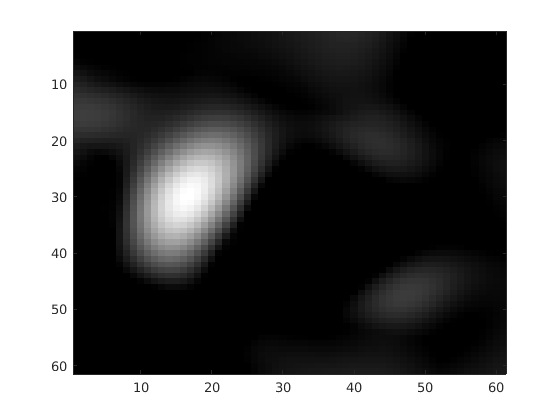}
  \caption{Regularized Kaczmarz}
\end{subfigure}\\
\begin{subfigure}{.23\textwidth}
  \centering
  \includegraphics[width=\linewidth]{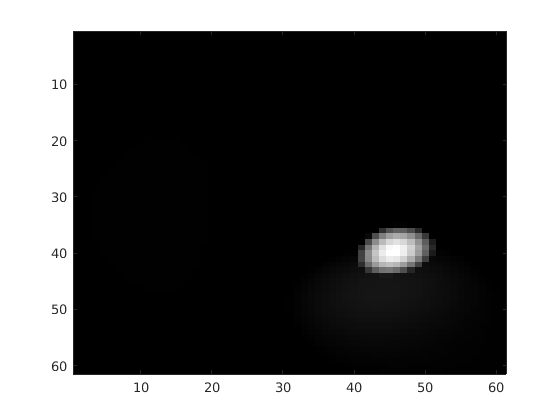}
  \caption{RESESOP-Kaczmarz}
\end{subfigure}%
\begin{subfigure}{0.23\textwidth}
  \centering
  \includegraphics[width=\linewidth]{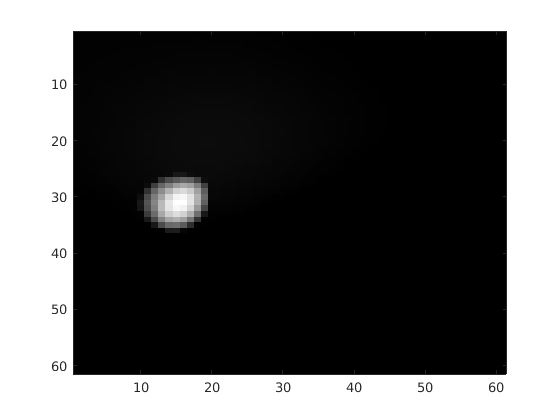}
  \caption{RESESOP-Kaczmarz}
\end{subfigure}\\
\caption{Dynamic Reconstructions of frame four from noisy simulated data with two different rotation speeds. The left column shows results for a rotation speed of 44 fpr, the right column for seven fpr.}
\label{fig:rec_sim}
\end{figure}

The reconstruction results are illustrated in \autoref{fig:rec_sim}.
The first row depicts the ground truth at frame four for both rotational velocities.
The second row shows the images reconstructed by regularized Kaczmarz. 
In both dynamic scenarios significant artifacts are visible. Since the algorithm cannot rely on averaging over multiple frames, the noise in the data is not satisfactorily reduced, resulting in artifacts throughout the entire field-of-view. 
Furthermore, the circular shape of the object is blurred, distorted and hence no more recognizable. As to be expected, this effect is worse the faster the motion. 
These examples motivate again that the dynamic behaviour of the concentration needs to be taken into account within the reconstruction step. 

The respective results of RESESOP-Kaczmarz are presented in the third row of figure \ref{fig:rec_sim}.
In both dynamic scenarios, the actual shape of the object is recognizable in the reconstructions, i.e. the algorithm does indeed compensate for the motion. However, due to the strict restriction to frame-sized subproblems, we can still observe some small distortions in case of the very fast motion. Also the location of the object is not correctly reconstructed but depicted in the middle of the traveled path during frame four and not in the beginning. However, we further note that in both dynamic cases, the method is able to eliminate the artifacts caused by the noise throughout the field of view.

Altogether, \autoref{fig:rec_sim} demonstrates the advantage of applying RESESOP-Kaczmarz to dynamic MPI problems. 

\subsection{Parameters for RESESOP-Kaczmarz}
In this section, we want to study the robustness of the RESESOP-Kaczmarz algorithm regarding estimates of the total inexactness levels $\zeta_{i}$, number of iterations and size of subproblems. 
All experiments in this subsection have been conducted with simulated data for the fast rotating object with ten percent of added noise.
The reconstruction in \autoref{fig:rec_sim} (f) was obtained by determining $\zeta_{i}$ as described in Section II.III.1 with the Euclidean norm and by performing ten full RESESOP-Kaczmarz iterations while the size of all subproblems corresponds to one full frame.

To test stability regarding the estimates of the inexactness levels, we performed the reconstruction with scaled versions of these inexactness levels (more precisely with factors $10^{-1}, \, 10^{1}$ and $10^5$ respectively) representing an under-, respectively over-estimations. In addition, we also considered 15\% noise added to the norm estimates $\zeta_i$ to test robustness with respect to computation errors. Lastly, we further use the alternative approach based on prior reconstruction discussed in section \ref{sec:in}.

\begin{table}[ht]
\setlength{\tabcolsep}{-1pt}
\centering
\begin{tabular}{c@{\hskip -2pt}ccc}
& Norm & 15\% Noise & RegKacz\\
$\zeta \cdot 10^{-1}$&\includegraphics[scale=.165,align=c]{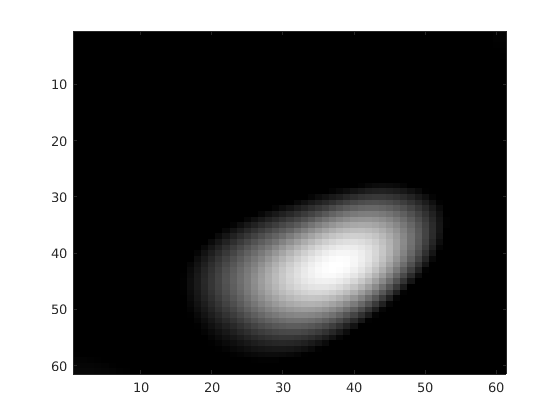}&\includegraphics[scale=.165,align=c]{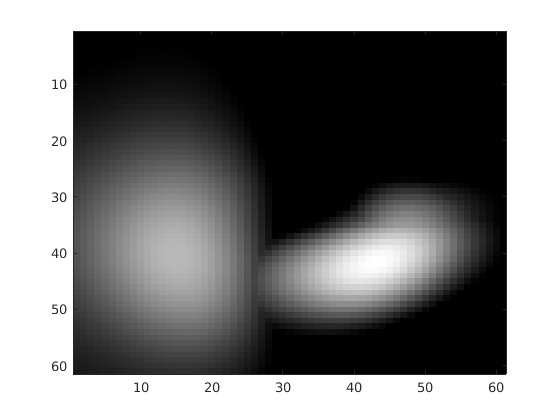}&\includegraphics[scale=.165,align=c]{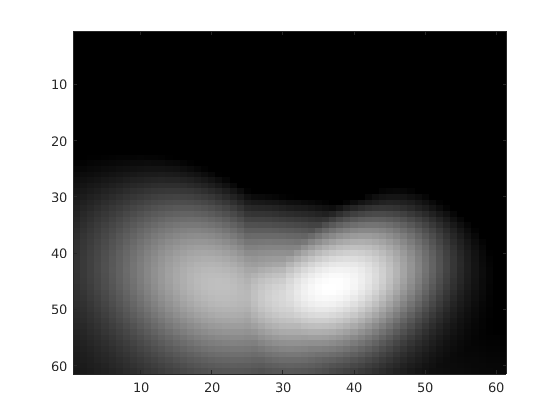}\\
$\zeta \cdot 10^{0}$&\includegraphics[scale=.165,align=c]{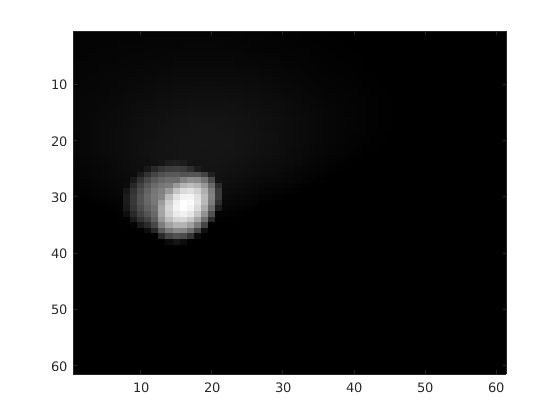}&\includegraphics[scale=.165,align=c]{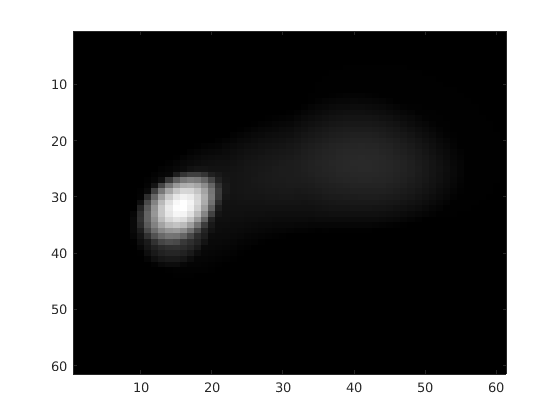}&\includegraphics[scale=.165,align=c]{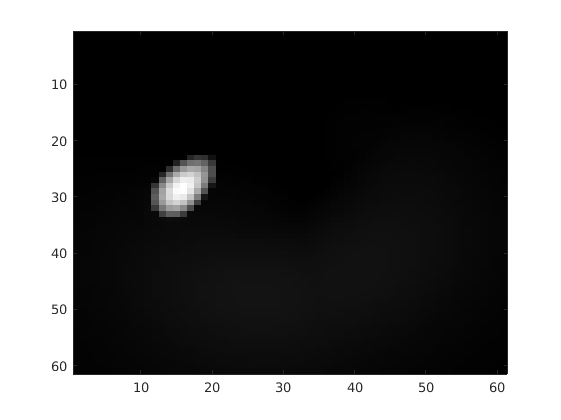}\\
$\zeta \cdot 10^{1}$&\includegraphics[scale=.165,align=c]{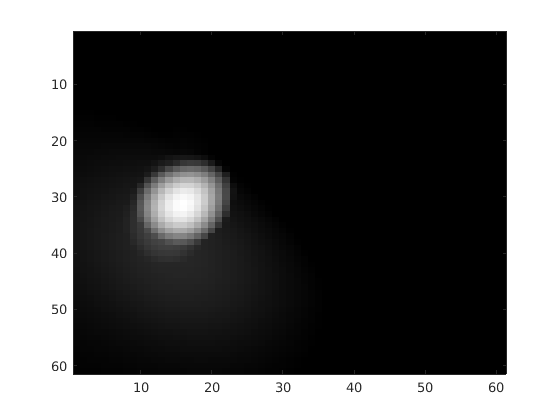}&\includegraphics[scale=.165,align=c]{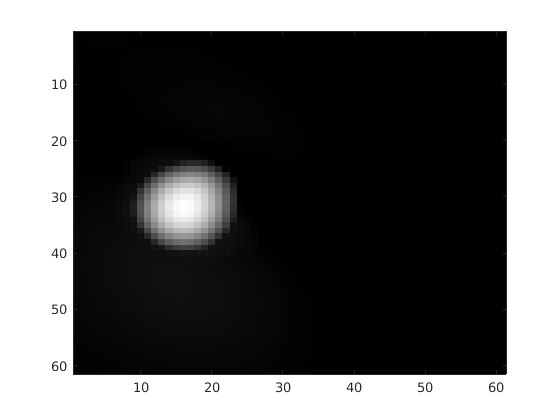}&\includegraphics[scale=.165,align=c]{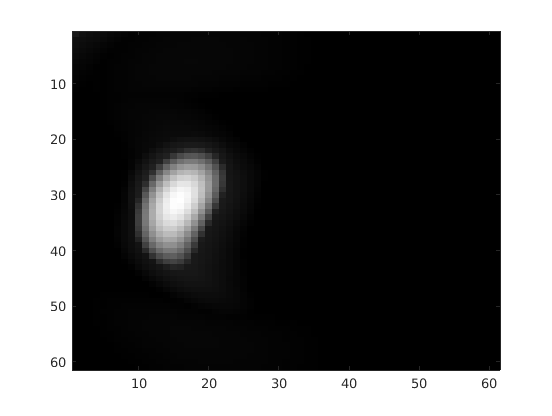}\\
$\zeta \cdot 10^{5}$&\includegraphics[scale=.165,align=c]{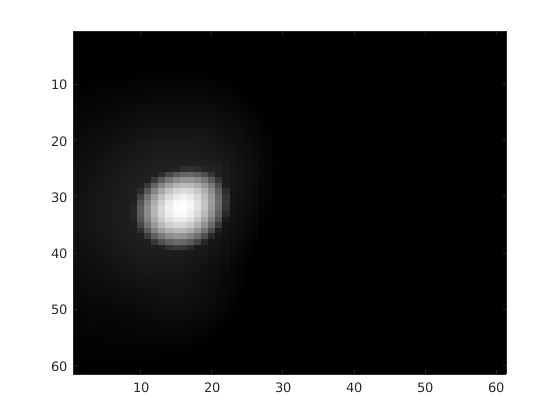}&\includegraphics[scale=.165,align=c]{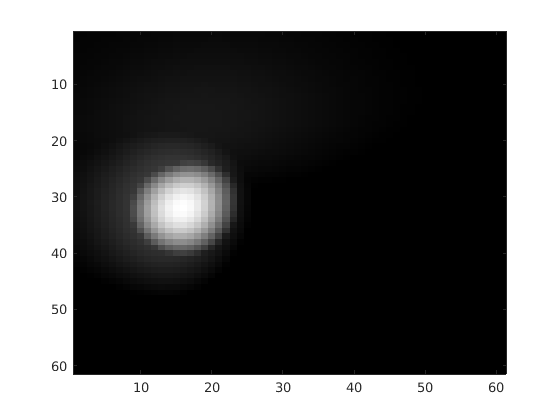}&\includegraphics[scale=.165,align=c]{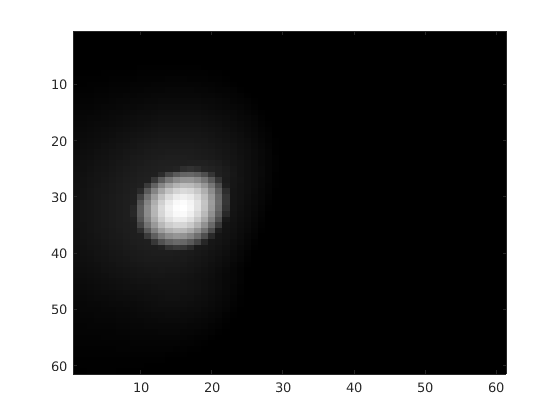}\\
\end{tabular}
\caption{Dynamic Reconstructions with RESESOP-Kaczmarz for different estimates of the total inexactness levels $\zeta_i$. Column 1:
computed as distance of two frames via euclidean norm, Column 2 : as in Column 1 with 15 \% added noise , Column 3: computed as MSE from regularized Kaczmarz reconstructions. 
Each row corresponds to one additional scaling factor.}
\label{tab:uncertainties}
\end{table}

\autoref{tab:uncertainties} illustrates the influence of different inexactness levels on the reconstruction result of RESESOP-Kaczmarz.
Both the proposed norm estimation directly from the measured data as well as the estimation from prior static reconstructions with regularized Kaczmarz result in images of similar quality (unless in the underestimated case, see first row, where the computation from the data is beneficial).
Nevertheless, an underestimation of the inexactness levels leads to an increase in motion artefacts, while the motion is well compensated for even if largely overestimated levels are used. 
The algorithm is also stable regarding (random) errors in the computed inexactness levels, see Column 2.
Despite strong deviations of overall 15\%, the algorithm still provides a reconstruction of comparable image quality as with the exact values.

\begin{figure}[!htb]
	\centering
	\includegraphics[scale=0.5]{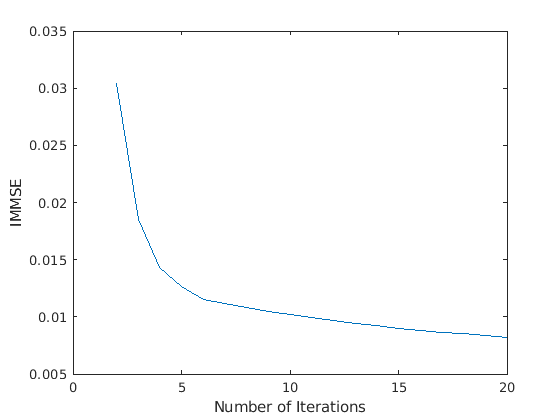}
	\caption{Mean squared error in dependence on the amount of RESESOP-Kaczmarz loops.}
	\label{fig:Iterations}
\end{figure}

Next, we study how the reconstruction quality develops in dependence on the number of iterations. 
As was shown in \cite{Blanke2020}, the algorithm converges monotone as long as adequate search directions are chosen.
This is clearly visible in figure \ref{fig:Iterations} which shows the evolution of the mean-squared error of the RESESOP-Kaczmarz solution compared to the ground truth with increasing number of iterations.
Furthermore, we observe that the error reduces rapidly during the first iteration and then flattens out.

\begin{figure}[!htb]
\centering
\begin{subfigure}{.16\textwidth}
  \centering
  \includegraphics[width=\linewidth]{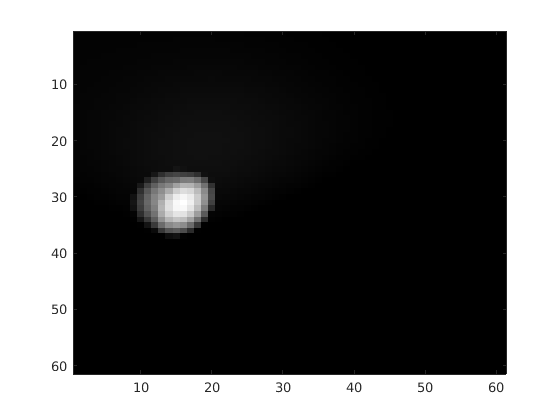}
  \caption{Frame}
\end{subfigure}%
\hspace{-.2cm}
\begin{subfigure}{0.16\textwidth}
  \centering
  \includegraphics[width=\linewidth]{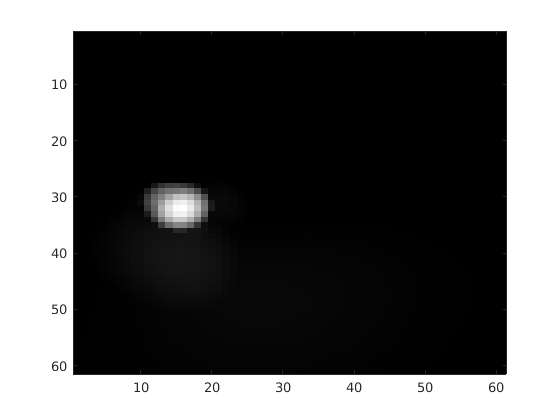}
  \caption{$\frac{1}{2}$-frame}
  \label{fig:rec_sim_fast_sub}
\end{subfigure}
\hspace{-.2cm}
\begin{subfigure}{.16\textwidth}
  \centering
  \includegraphics[width=\linewidth]{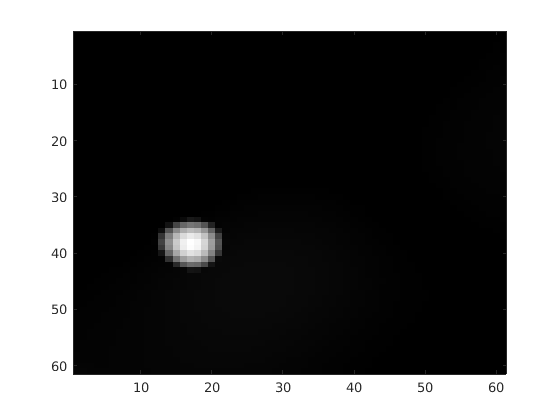}
  \caption{$\frac{1}{4}$-frame}
  \label{fig:quarter}
\end{subfigure}\\
\begin{subfigure}{0.16\textwidth}
  \centering
  \includegraphics[width=\linewidth]{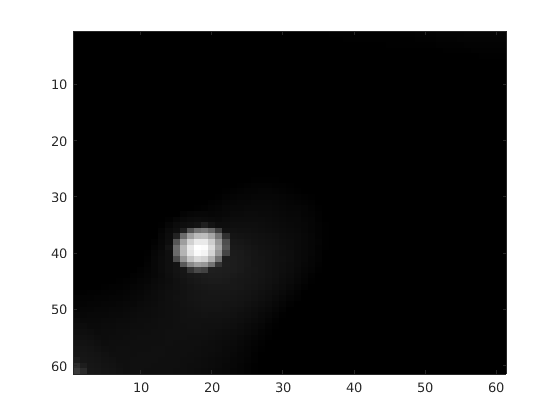}
  \caption{$\frac{1}{8}$-frame}
\end{subfigure}
\hspace{-.2cm}
\begin{subfigure}{.16\textwidth}
  \centering
  \includegraphics[width=\linewidth]{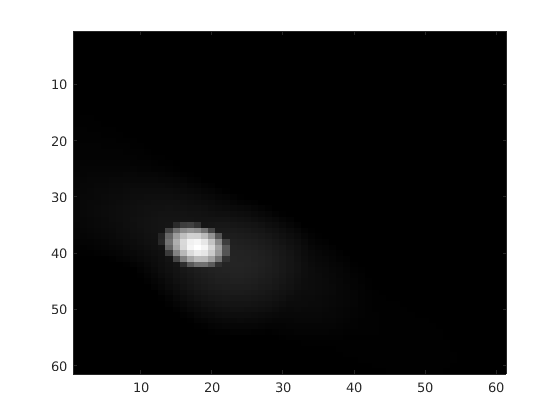}
  \caption{$\frac{1}{16}$-frame}
  \label{fig:16}
\end{subfigure}%
\hspace{-.2cm}
\begin{subfigure}{0.16\textwidth}
  \centering
  \includegraphics[width=\linewidth]{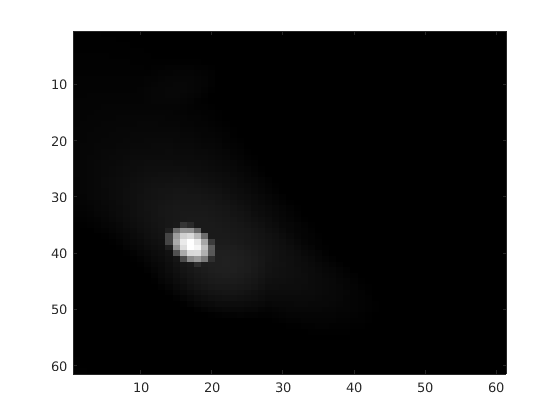}
  \caption{$\frac{1}{32}$-frame}
  \label{fig:32}
\end{subfigure}
\caption{Dynamic Reconstructions from noisy simulated data in the fast scenario where different sizes of subproblems are considered within the RESESOP-Kaczmarz algorithm.
}
\label{fig:rec_sim_sub}
\end{figure}

Lastly, we want to examine the influence of the size of the subproblems.
Depending on the speed of the motion, assuming the concentration to be static during a complete Lissajous trajectory can be too much simplified.
Indeed, \autoref{fig:rec_sim} shows that under this assumption in the very fast scenario, the reconstructed object is not correctly located.
This can be improved by considering subproblems of smaller size, see \autoref{fig:rec_sim_sub}.

Choosing the sizes of the subproblems as a quarter of a frame or smaller improves the reconstructed location of the object, see e.g. \autoref{fig:quarter}.
However, choosing the size too small can reintroduce a low amount of noise in the reconstruction, see  Figures \ref{fig:16} and \ref{fig:32}.
Thus, in case of very fast motions, one has to balance motion compensation and noise reduction when choosing appropriate sizes of subproblems. 

\subsection{Real Data}
In this section we evaluate the performance of RESESOP-Kaczmarz on the real data introduced in Section \ref{RealData}.

\begin{table}[ht]
\setlength{\tabcolsep}{-1pt}
\centering
\begin{tabular}{c@{\hskip 2pt}cc}
& Frame & $\frac{1}{4}$-frame \\
Frame 1&\includegraphics[scale=.22,align=c]{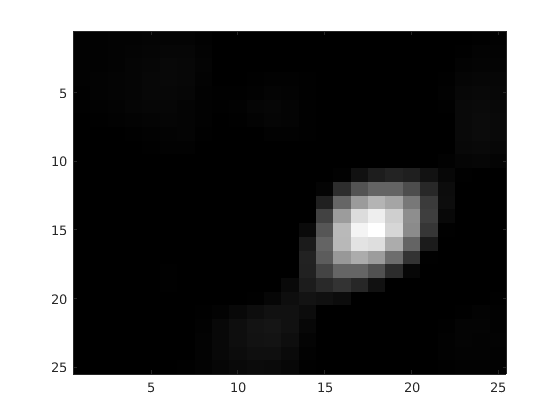}&\includegraphics[scale=.22,align=c]{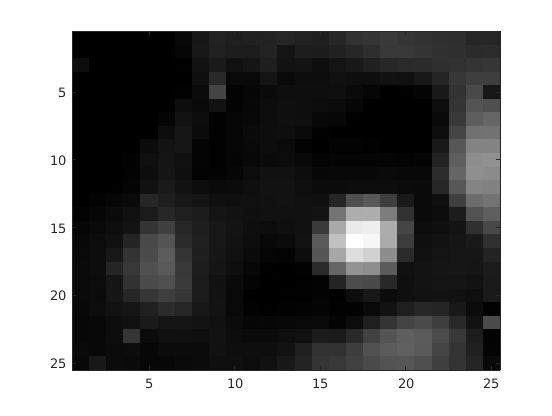}\\
Frame 2&\includegraphics[scale=.22,align=c]{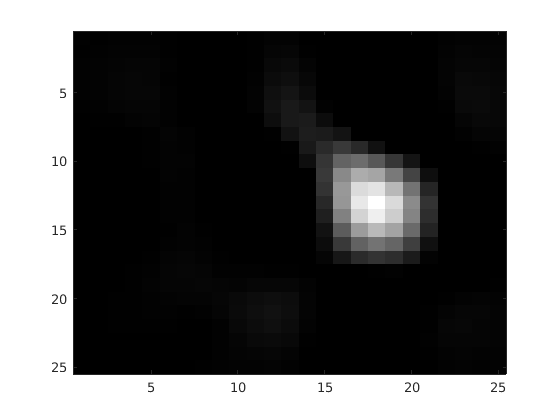}&\includegraphics[scale=.22,align=c]{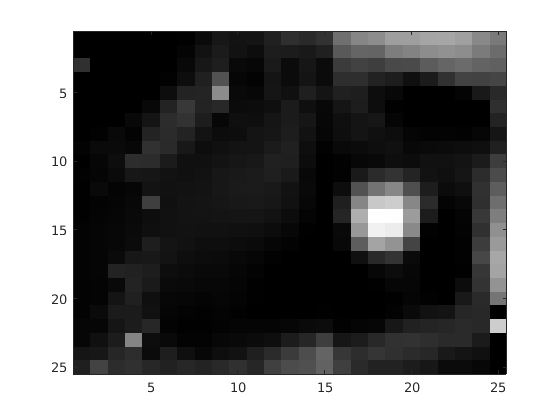}\\
Frame 3&\includegraphics[scale=.22,align=c]{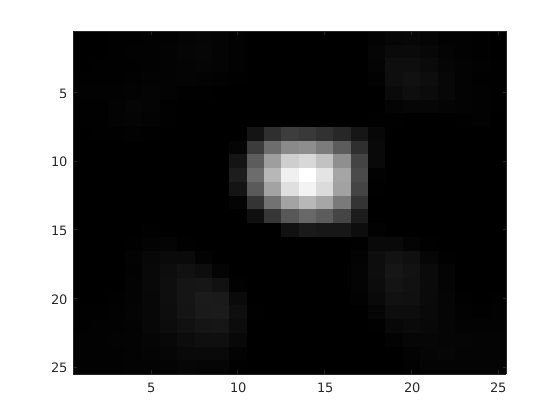}&\includegraphics[scale=.22,align=c]{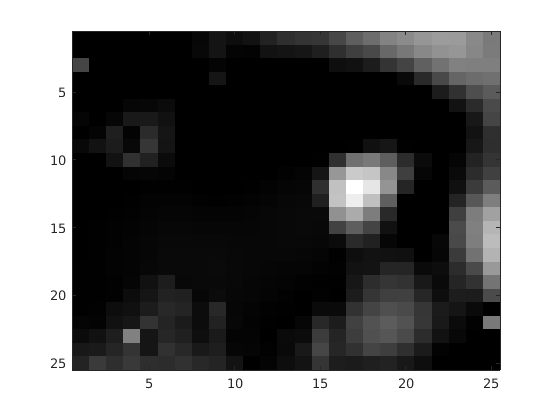}\\
Frame 4&\includegraphics[scale=.22,align=c]{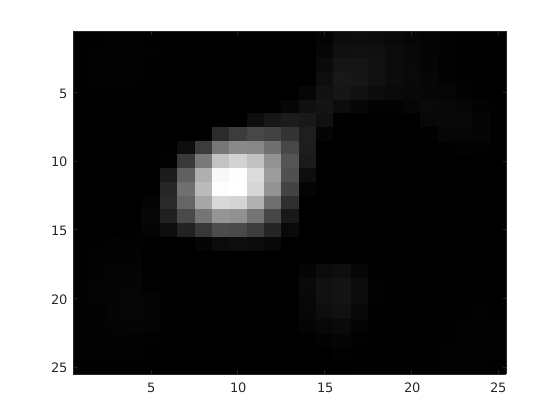}&\includegraphics[scale=.22,align=c]{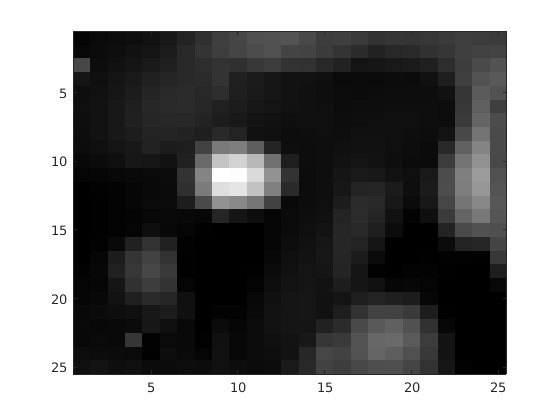}
\end{tabular}
\caption{Dynamic Reconstructions with RESESOP-Kaczmarz of four consecutive frames of real data for a 7 Hz rotation with different sized subproblems. The original numeration of the frames starts at frame 750.}
\label{tab:real}
\end{table}

\autoref{tab:real} depicts reconstructions of four consecutive frames of real data.
The left column shows images computed with RESESOP-Kaczmarz using frame-sized subproblems. From the reconstructed images, we can clearly deduce the temporal evolution of the object, namely the rotation of the glass capillary. Also its spatial location per time step is overall well captured and the noise well suppressed. However, examining frames one, two and four the circular shape of the reconstructed glass capillary is slightly distorted due to the fast motion. 

The right column of \autoref{tab:real} shows the respective reconstruction results using $\frac{1}{4}$-frames as subproblems. Compared to the frame-wise case, they portray the object with a more constant round shape but are also more affected by noise. This illustrates again the trade-off between noise reduction and motion compensation. In a future step, it would be interesting to study whether algorithms (or reconstructions) for different subproblem sizes could be merged in order to combine the good motion compensation property of small subproblems with the better noise reduction of large subproblems.  

\section{Conclusion}
In this article, we propose the RESESOP-Kaczmarz algorithm to reconstruct dynamic objects from MPI data. 
The method takes motion into account as a model inexactness and, within the MPI framework, the only required a priori information can be computed directly from the measured data. The potential of the method was demonstrated on real and simulated data. In particular, our detailed experiments on simulated data show the robustness of the method regarding its various parameters and that it outperforms the regularized Kaczmarz algorithm, a common solver in MPI based on a stationary assumption on the object. Since the method further allows to consider subproblems smaller than one complete frame, it is suitable for a variety of dynamic MPI problems, including scenarios with rapid particle movements. 

\section*{Acknowledgments}
The authors acknowledge the support by 'Deutsche Forschungsgemeinschaft' under grant HA 8176/2-1. The authors further thank Hannes Albers and Tobias Kluth for providing us with simulated data as well as Christina Brandt and Tobias Knopp for real data.

\printbibliography

\end{document}